\documentclass[12pt]{amsart}
\usepackage{fullpage,url}
\usepackage[dvips]{hyperref}
 \usepackage{amscd}   
\usepackage{amsmath}

\DeclareFontEncoding{OT2}{}{} 

\DeclareMathOperator{\per}{per}

\DeclareMathOperator{\tor}{tor}
\DeclareMathOperator{\alg}{alg}
\DeclareMathOperator{\hhat}{\hat{h}}

\DeclareMathOperator{\rk}{rk}

\newcommand{\tensor}{\otimes}

\newtheorem{theorem}{Theorem}[section]
\newtheorem{lemma}[theorem]{Lemma}
\newtheorem{corollary}[theorem]{Corollary}

\theoremstyle{definition}
\newtheorem{definition}[theorem]{Definition}

\theoremstyle{remark}
\newtheorem{remark}[theorem]{Remark}

\newtheorem{sublemma}[theorem]{Sublemma}

\title{Elliptic curves over the perfect closure of a function field}
\author{Dragos Ghioca}
\begin{document}
\begin{abstract}
We prove that the group of rational points of a non-isotrivial elliptic curve defined over the perfect closure of a function field in one variable over a finite field is finiteley generated.
\end{abstract}
\maketitle

\section{Introduction}

The main ingredient for our result is a study of the Lehmer inequality for elliptic curves. The classical Lehmer conjecture (see \cite{Leh}, page $476$) asserts that there is
an absolute constant $C>0$ 
so that any algebraic number $\alpha$ that is not a
root of unity satisfies the following inequality for its logarithmic height
$$h(\alpha)\ge\frac{C}{[\mathbb{Q}(\alpha):\mathbb{Q}]}.$$
A partial result 
towards this conjecture is obtained in \cite{Dob}. The
analog of Lehmer 
conjecture for elliptic curves and abelian varieties asks for a good lower bound for the canonical height of a non-torsion point of the abelian variety. Also this question has been much studied 
(see \cite{Siabba},
\cite{DaHi}, \cite{HiSi}, \cite{Mass}, \cite{Siabell}). As a consequence of our work, we obtain the following result.

\begin{theorem}
\label{T:ellper}
Let $K$ be a function field of transcendence degree
$1$ over $\mathbb{F}_p$. Let 
$E$ be a non-isotrivial elliptic curve defined over
$K$. Then $E(K^{\per})$ is 
finitely generated.
\end{theorem}

Using completely different methods, Minhyong Kim studied the set of rational points of non-isotrivial curves of genus at least two over the perfect closure of a function field in one variable over a finite field (see \cite{kim}).

Using the result of Theorem~\ref{T:ellper}, Thomas
Scanlon proved the full 
positive characteristic Mordell-Lang Conjecture for
abelian varieties that are 
isogenous to a product of elliptic curves (see \cite{Sca2}).

Our plan is to prove first some general results about tame modules in Section~\ref{se:tame} and then to prove Theorem~\ref{T:ellper} in Section~\ref{se:ellper}.
 
\section{Tame modules}
\label{se:tame}

\begin{definition}
\label{D:tamemodules}
Let $R$ be an integral domain and let $K$ be its field
of fractions. If $M$ is an $R$-module, then by the
\emph{rank} of $M$, denoted $\rk(M)$, we mean the dimension
of the $K$-vector space $M\tensor_{R}K$. We call $M$ a
\emph{tame} module if every finite rank submodule of
$M$ is finitely generated.
\end{definition}

\begin{lemma}
\label{L:heightonmodules}
Let $R$ be a Dedekind domain and let $M$ be an
$R$-module. Assume there exists a function
$h:M\rightarrow\mathbb{R}_{\ge 0}$ satisfying the
following properties

(i) (triangle inequality) $h(x\pm y)\le h(x)+h(y)$,
for every $x,y\in M$.

(ii) if $x\in M_{\tor}$, then $h(x)=0$.

(iii) there exists $c>0$ such that for each $x\notin
M_{\tor}$, $h(x)>c$.

(iv) there exists $a\in R\setminus\{0\}$ such that $R/aR$ is finite and for all 
$x\in M$, $h(ax)\ge 4h(x)$.

If $M_{\tor}$ is finite, then $M$ is a tame
$R$-module.
\end{lemma}

\begin{proof}
By the definition of a tame module, it suffices to
assume that $M$ is a finite rank $R$-module and
conclude that it is finitely generated.

Let $a\in R$ as in $(iv)$ of Lemma
\ref{L:heightonmodules}. By Lemma $3$ of \cite{Poo},
$M/aM$ is finite (here we use the assumption that
$M_{\tor}$ is finite). The following result is the key
to the proof of Lemma~\ref{L:heightonmodules}.

\begin{sublemma}
\label{S:finitely}
For every $D>0$, there exists finitely many $x\in M$
such that $h(x)\le D$.
\end{sublemma}

\begin{proof}[Proof of Sublemma \ref{S:finitely}.]

If we suppose Sublemma~\ref{S:finitely} is not true,
then we can define 
$$C=\inf\{D\mid\text{ there exists infinitely many
$x\in M$ such that 
$h(x)\le D$}\}.$$
Properties $(ii)$ and $(iii)$ and the finiteness of
$M_{\tor}$ yield $C\ge c>0$. By the definition of $C$,
it must be that there exists an infinite sequence of
elements $z_n$ of $M$ such that for every $n$, 
$$h(z_n)<\frac{3C}{2}.$$

Because $M/aM$ is finite, there exists a coset of $aM$
in $M$ containing infinitely many $z_n$ from the above
sequence.

But if $k_1\ne k_2$ and $z_{k_1}$ and $z_{k_2}$ are in
the same coset of $aM$ in $M$, then let 
$y\in M$ be such that $ay=z_{k_1}-z_{k_2}$. Using
properties $(iv)$ and $(i)$, we get
$$h(y)\le\frac{h(z_{k_1}-z_{k_2})}{4}\le\frac{h(z_{k_1})+h(z_{k_2})}{4}<\frac{3C
}{4}.$$
We can do this for any two elements of the sequence
that lie in the same 
coset of $aM$ in $M$. Because there are infinitely
many of them lying in 
the same coset, we can construct infinitely many
elements $z\in M$ such that $h(z)<\frac{3C}{4}$,
contradicting the minimality of $C$.
\end{proof}

From this point on, our proof of Lemma \ref{L:heightonmodules} follows the 
classical
descent argument in the 
Mordell-Weil theorem (see \cite{Ser}).

Take coset representatives $y_1,\dots,y_k$ for $aM$ in
$M$. Define then 
$$B=\max_{i\in\{1,\dots,k\}}h(y_i).$$
Consider the set $Z=\{x\in M\mid h(x)\le B\}$, which
is finite 
according to Sublemma~\ref{S:finitely}. Let 
$N$ be the finitely generated $R$-submodule of $M$
which is spanned by $Z$.

We claim that $M=N$. If we suppose this is not
the case, then by Sublemma~\ref{S:finitely} we can
pick $y\in M-N$ which minimizes $h(y)$. Because $N$
contains all the coset representatives of 
$aM$ in $M$, 
we can find $i\in\{1,\dots,k\}$ such that $y-y_i\in
aM$. Let $x\in M$ be 
such that
$y-y_i=ax$. Then $x\notin N$ because 
otherwise it would follow that $y\in N$ ( we already know $y_i\in N$).
By our choice of $y$ and by properties $(iv)$ and
$(i)$, we have
$$h(y)\le
h(x)\le\frac{h(y-y_i)}{4}\le\frac{h(y)+h(y_i)}{4}
\le\frac{h(y)+B}{4}.$$
This means that $h(y)\le\frac{B}{3}<B$. This
contradicts the fact that 
$y\notin N$ because $N$ contains all the elements
$z\in M$ such that $h(z)\le B$. This contradiction
shows that indeed $M=N$ and so, $M$ is finitely 
generated. 
\end{proof}

\begin{corollary}
\label{C:alephcardinal}
Let $R$ be a Dedekind domain and let $M$ be a tame
$R$-module. 

(a) If $\rk(M)=\aleph_0$, then $M$ is a direct sum of
a finite torsion submodule and a free submodule of rank
$\aleph_0$.

(b) If $\rk(M)$ is finite and $R$ is a principal ideal
domain, then $M$ is a direct sum of a finite torsion
submodule and a free submodule of finite rank.
\end{corollary}

\begin{proof}
Part $(a)$ of Corollary \ref{C:alephcardinal} is
proved in Proposition $10$ of \cite{Poo}.

If $\rk(M)$ is finite and because $M$ is tame, we
conclude that $M$ is finitely generated. Because $R$
is a principal ideal domain we get the result of part
$(b)$ of Corollary ~\ref{C:alephcardinal}.
\end{proof}

\section{Elliptic curves}
\label{se:ellper}

The setting is the following: $K$ is
a finitely generated 
field of transcendence degree $1$ over $\mathbb{F}_p$
where $p$ is a prime as always. We fix an algebraic
closure $K^{\alg}$ of $K$. We denote by
$\mathbb{F}_p^{\alg}$ the algebraic closure of
$\mathbb{F}_p$ inside $K^{\alg}$.

Let $E$ be a non-isotrivial elliptic curve (i.e. 
$j(E)\notin\mathbb{F}_p^{\alg}$) defined over $K$. Let
$K^{\per}$ be the 
perfect closure of $K$ inside $K^{\alg}$. We will
prove in Theorem~\ref{T:ellper} that $E(K^{\per})$ is
finitely generated. 

For every finite extension $L$ of $K$ we denote by
$M_L$ the set of discrete valuations $v$ on $L$,
normalized so that the value group of $v$ is
$\mathbb{Z}$. For each $v\in M_L$ we denote by $f_v$
the degree of the residue field of $v$ over
$\mathbb{F}_p$. If $P\in E(L)$ and $m\in\mathbb{Z}$,
$mP$ represents the point on the elliptic curve obtained using the group law on 
$E$. We define a notion of height for the point $P\in 
E(L)$
with respect to the field $K$ (see Chapter $VIII$ of \cite{Sil0} and Chapter 
$III$ of \cite{Sil})
\begin{equation}
\label{E:defheightelliptic}
h_{K}(P)=\frac{-1}{[L:K]}\sum_{v\in
M_L}f_v\min\{0,v(x(P))\}.
\end{equation}
Then we define the canonical height of $P$ with
respect to $K$ as
\begin{equation}
\label{E:defcanonicheightelliptic}
\hhat_{E/K}(P)=\frac{1}{2}\lim_{n\rightarrow\infty}\frac{h_K(2^nP)}{4^n}.
\end{equation}

We also denote by $\Delta_{E/K}$ the divisor which is the minimal
discriminant of $E$ with respect to the field $K$ (see
Chapter $VIII$ of \cite{Sil0}). By $\deg(\Delta_{E/K})$
we denote the degree of the divisor $\Delta_{E/K}$ (computed with respect to 
$\mathbb{F}_p$). We
denote by $g_K$ the genus of the function field $K$.

The following result is proved in \cite{GSz}.
\begin{theorem}[Goldfeld-Szpiro]
\label{T:theorem 7}
Let $E$ be an elliptic curve over a function field $K$ of
one variable over a field in any characteristic. Let
$\hhat_{E/K}$ denote the canonical height on $E$ and
let $\Delta_{E/K}$ be the minimal discriminant of $E$,
both computed with respect to $K$. Then for every
point $P\in E(K)$ which is not a torsion point:
$$\hhat_{E/K}(P)\ge 10^{-13}\deg(\Delta_{E/K})\text{
if }\deg(\Delta_{E/K})\ge 24(g_K-1)\text{ ,}$$
and
$$\hhat_{E/K}(P)\ge
10^{-13-23g}\deg(\Delta_{E/K})\text{ if
}\deg(\Delta_{E/K})<24(g_K-1).$$
\end{theorem}

We are ready to prove our result.
\begin{proof}[Proof of Theorem~\ref{T:ellper}.]  
We first observe that replacing $K$ by a finite extension does not affect the 
conclusion of the theorem. Thus, at the expense of replacing $K$ by a finite 
extension, we may assume $E$ is semi-stable over $K$ (the existence of such a 
finite extension is guaranteed by Proposition $5.4$ $(a)$ of \cite{Sil0}).

As before, we let $\hhat_{E/K}$ and $\Delta_{E/K}$ be
the canonical height on $E$ and the minimal
discriminant of $E$, respectively, computed with
respect to $K$.

For every $n\ge 1$, we denote by $E^{(p^n)}$ the
elliptic curve $F^n(E)$, where $F$ is the usual
Frobenius (seen as morphism of varieties). Thus 
\begin{equation}
\label{E:equation0}
F^n:E(K^{1/p^n})\rightarrow E^{(p^n)}(K)
\end{equation}
is a bijection. Moreover, for every $P\in
E\left(K^{1/p}\right)$, 
\begin{equation}
\label{E:ellipticconnection}
pP=\left(V F\right)\left(P\right) \in
V\left(E^{\left(p\right)}\left(K\right)\right)\subset
E\left(K\right)
\end{equation}
where $V$ is the Verschiebung. Similarly, we get that 
\begin{equation}
\label{E:elliptictower}
p^nE\left(K^{1/p^n}\right)\subset
E\left(K\right)\text{ for every $n\ge 1$.}
\end{equation}
Thus $E(K^{\per})$ lies in the $p$-division hull of
the $\mathbb{Z}$-module $E(K)$. Because $E(K)$ is
finitely generated (by the Mordell-Weil theorem), we
conclude that $E(K^{\per})$, as a $\mathbb{Z}$-module,
has finite rank. 

We will show next that the height function
$\hhat_{E/K}$ and $p\in\mathbb{Z}$ satisfy the
properties $(i)$-$(iv)$ of Lemma
~\ref{L:heightonmodules}. Properties $(i)$ and $(ii)$
are well-known for $\hhat_{E/K}$ and we also have the
formula (see Chapter $VIII$ of \cite{Sil0})
$$\hhat_{E/K}(pP)=p^2\hhat_{E/K}(P)\ge 4\hhat(P)\text{
for every $P\in E(K^{\alg})$,}$$
which proves that property $(iv)$ of Lemma~\ref{L:heightonmodules} holds.

Let $P$ be a non-torsion point of $E(K^{\per})$. Then
$P\in E(K^{1/p^n})$ for some $n\ge 0$. Because
$K^{1/p^n}$ is isomorphic to $K$, they have the same
genus, which we call it $g$. We denote by
$\hhat_{E/K^{1/p^n}}$ and $\Delta_{E/K^{1/p^n}}$ the
canonical height on $E$ and the minimal discriminant
of $E$, respectively, computed with respect to
$K^{1/p^n}$. Using Theorem~\ref{T:theorem 7}, we
conclude
\begin{equation}
\label{E:heightineq7_n}
\hhat_{E/K^{1/p^n}}(P)\ge
10^{-13-23g}\deg(\Delta_{E/K^{1/p^n}}).
\end{equation}

We have
$\hhat_{E/K^{1/p^n}}(P)=[K^{1/p^n}:K]\hhat_{E/K}(P)=p^n\hhat_{E/K}(P)$.
Similarly, using the proof of Proposition $5.4$ $(b)$ of \cite{Sil0}, 
$$\deg(\Delta_{E/K^{1/p^n}})=p^n\deg(\Delta_{E/K}).$$ 
We
conclude that for every $P\in E(K^{\per})$,
\begin{equation}
\label{E:heightineq7}
\hhat_{E/K}(P)\ge 10^{-13-23g}\deg(\Delta_{E/K}).
\end{equation}
Because $E$ is non-isotrivial, $\Delta_{E/K}\ne 0$ and
so, $\deg(\Delta_{E/K})\ge 1$. We conclude
\begin{equation}
\label{E:heightineq7'}
\hhat_{E/K}(P)\ge 10^{-13-23g}.
\end{equation}
Inequality \eqref{E:heightineq7'} shows that property
$(iii)$ of Lemma \ref{L:heightonmodules} holds for
$\hhat_{E/K}$. Thus properties $(i)$-$(iv)$ of Lemma
~\ref{L:heightonmodules} hold for $\hhat_{E/K}$ and
$p\in\mathbb{Z}$.

We show that $E_{\tor}(K^{\per})$ is finite. Equation
\eqref{E:elliptictower} shows that the
prime-to-$p$-torsion of $E(K^{\per})$ equals the
prime-to-$p$-torsion of $E(K)$; thus the
prime-to-$p$-torsion of $E(K^{\per})$ is finite. If
there exists
infinite $p$-power torsion in
$E\left(K^{\per}\right)$, equation \eqref{E:equation0}
yields that we have 
arbitrarily large $p$-power torsion in the family of
elliptic curves 
$E^{\left(p^n\right)}$ over $K$. But this contradicts
standard results on 
uniform boundedness for the torsion of elliptic curves
over function fields, as 
established in \cite{Lev} (actually, \cite{Lev} proves a uniform boundedness of 
the entire torsion of elliptic curves over a fixed function field; thus 
including the prime-to-$p$-torsion). Hence $E_{\tor}(K^{\per})$
is finite.

Because all the hypothesis of Lemma~\ref{L:heightonmodules} hold, we conclude that
$E(K^{\per})$ is tame. Because
$\rk\left(E(K^{\per})\right)$ is finite we conclude by
Corollary ~\ref{C:alephcardinal} $(b)$ that
$E(K^{\per})$ is a direct sum of a finite torsion
submodule and a free submodule of finite rank.
\end{proof}

\begin{remark}
It is absolutely crucial in Theorem \ref{T:ellper} that
$E$ is non-isotrivial. Theorem~\ref{T:ellper} fails in
the isotrivial case, i.e. there exists no $n\ge 0$ such that 
$E(K^{\per})=E(K^{1/p^n})$. Indeed, if $E$ is defined by
$y^2=x^3+x$ ($p>2$),
$K=\mathbb{F}_p\left(t,(t^3+t)^{\frac{1}{2}}\right)$
and $P=\left(t,(t^3+t)^{\frac{1}{2}}\right)$, then
$F^{-n}P\in 
E(K^{1/p^n})\setminus E(K^{1/p^{n-1}})$, for every $n\ge 1$. So, $E(K^{\per})$
is not 
finitely generated in this case (and we can get a similar
example also for the case $p=2$).
\end{remark}

\end{document}